\documentclass[12pt,a4paper,reqno]{amsart}

\usepackage[utf8]{inputenc}
\usepackage{enumerate,bbm,titletoc,amssymb,amsthm,mathabx,graphicx,mathrsfs}
\usepackage{xspace}
\usepackage[normalem]{ulem}
\usepackage{euflag}

\newcommand{\ratQM}{rational quadratic module\xspace}

\newcommand{\cH}{\mathcal{H}}
\newcommand{\cK}{\mathcal{K}}
\newcommand{\cE}{\mathcal{E}}
\newcommand{\cS}{\mathcal{S}}
\newcommand{\bb}{\mathbbm}
\newcommand{\mF}{\mathbbm F}
\newcommand{\mG}{\mathbbm G}
\newcommand{\ee}{\mathbbm e}
\newcommand{\one}{\mathbbm 1}
\DeclareMathOperator{\UCP}{UCP}
\DeclareMathOperator{\CP}{CP}
\usepackage{tikz}
\usepackage[hidelinks]{hyperref}

\DeclareMathOperator{\id}{id}

\DeclareMathOperator{\GL}{GL}
\DeclareMathOperator{\dom}{dom}
\DeclareMathOperator{\hdom}{hdom}
\DeclareMathOperator{\spa}{span}
\DeclareMathOperator{\opm}{M}
\DeclareMathOperator{\oph}{H}

\newcommand{\ve}{\varepsilon}

\newcommand{\N}{\mathbb{N}}

\newcommand{\R}{\mathbb{R}}
\newcommand{\C}{\mathbb{C}}
\newcommand{\cA}{\mathcal{A}}
\newcommand{\cB}{\mathcal{B}}
\newcommand{\cD}{\mathcal{D}}

\newcommand{\cO}{\mathcal{O}}

\newcommand{\cR}{\mathcal{R}}

\newcommand{\cT}{\mathcal{T}}

\newcommand{\cW}{\mathcal{W}}

\newcommand{\rr}{\mathbbm r}
\newcommand{\rs}{\mathbbm s}

\newcommand{\qq}{\mathbbm q}
\newcommand{\uu}{\mathbbm u}
\newcommand{\rv}{\mathbbm v}
\newcommand{\fX}{\mathfrak{X}}

\newcommand*{\mathlarger}[1]{\scalebox{1.35}{$#1$}}
\newcommand*{\mat}[1]{\opm_{#1}(\C)}
\newcommand*{\gl}[1]{\GL_{#1}(\C)}
\newcommand*{\her}[1]{\oph_{#1}(\C)}
\newcommand{\cfm}{\opm_\infty(\C)}

\newcommand{\Langle}{\mathop{<}\!}
\newcommand{\Rangle}{\!\mathop{>}}

\newcommand{\px}{\C\!\Langle x\Rangle}

\makeatletter
\def\moverlay{\mathpalette\mov@rlay}
\def\mov@rlay#1#2{\leavevmode\vtop{
		\baselineskip\z@skip \lineskiplimit-\maxdimen
		\ialign{\hfil$#1##$\hfil\cr#2\crcr}}}
\makeatother

\newcommand{\plangle}{\moverlay{(\cr<}}
\newcommand{\prangle}{\moverlay{)\cr>}}
\newcommand{\rx}{\C\plangle x \prangle}
\newcommand{\rex}{\mathfrak{R}_{\C}(x)}

\def\gen{\texttt{g}}

\usepackage[margin=2.5cm]{geometry}

\newtheorem{thm}{Theorem}[section]
\newtheorem{lem}[thm]{Lemma}
\newtheorem{cor}[thm]{Corollary}
\newtheorem{prop}[thm]{Proposition}
\newtheorem{thmA}{Theorem}

\theoremstyle{definition}
\newtheorem{defn}[thm]{Definition}

\theoremstyle{remark}
\newtheorem{rem}[thm]{Remark}

\makeatletter
\@namedef{subjclassname@2020}{
	\textup{2020} Mathematics Subject Classification}
\makeatother

\numberwithin{equation}{section}

\title{Noncommutative Rational Positivstellensätze with Operator Coefficients}

\author[A.\ Jindal]{Abhay Jindal${}^{1,Q}$}
\address{Abhay Jindal, Faculty of Mathematics and Physics, University of Ljubljana, Slovenia}
\email{abhay.jindal@fmf.uni-lj.si}
\thanks{${}^1$Supported by the Slovenian Research Agency 
	program P1-0222 and grants J1-50002 and N1-0536.}

\author[I.\ Klep]{Igor Klep${}^{2,Q,R}$}
\address{Igor Klep, Faculty of Mathematics and Physics, University of Ljubljana 
	\& Famnit, University of Primorska, Koper 
	\& Institute of Mathematics, Physics and Mechanics,
	Ljubljana, Slovenia}
\email{igor.klep@fmf.uni-lj.si}
\thanks{${}^2$Supported by the Slovenian Research Agency 
	program P1-0222 and grants J1-50002,
	N1-0217,
	J1-60011, J1-50001, J1-3004 and J1-60025.}

\author[J.\ Vol\v ci\v c]{Jurij Vol\v ci\v c${}^{3}$}
\address{The University of Auckland, Department of Mathematics,
	Private Bag 92019,
	Auckland 1142,
	New Zealand.}
\email{jurij.volcic@auckland.ac.nz}
\thanks{${}^3$Supported by the Marsden Fund MFP-UOA2528 from Government funding, administered by the Royal Society Te Ap\=arangi, and the New Staff grant 3734441 by the Faculty of Science, the University of Auckland.}

\thanks{${}^Q$This work was performed within the project COMPUTE, funded within the QuantERA II 
	Programme that has received funding from the EU's H2020 research and innovation programme under the GA No.~101017733 {\normalsize\euflag}}

\thanks{${}^R$Part of this work was performed within the project SDPcode, funded within the QuantERA III 
	Programme that has received funding from the EU's Horizon research and innovation programme under the GA No.~101212998 {\normalsize\euflag}}



\begin{document}
	
	\begin{abstract}
		We establish a convex Positivstellensatz for noncommutative rational functions with bounded operator coefficients. A hermitian operator-valued rational function that is positive semidefinite on a matrix convex set described by a linear pencil admits a positivity certificate consisting of rational hermitian squares and a completely positive term involving the defining linear pencil. The certificate comes with an explicit complexity bound independent of the dimension of the coefficient Hilbert space. As a consequence, we solve the operator-valued noncommutative rational analog of Hilbert's 17th problem: every positive semidefinite operator-valued noncommutative rational function is a sum of hermitian squares, and even a single hermitian square when the coefficient Hilbert space is infinite-dimensional. In contrast, the corresponding commutative statement fails for infinite-dimensional operator coefficients.
	\end{abstract}
	
	\maketitle

	\tableofcontents

	\section{Introduction}\label{s:intro}
	
	Hilbert's 17th problem asks whether every nonnegative real polynomial is a
	sum of squares of rational functions. Artin's affirmative solution
	is a cornerstone of real algebraic geometry \cite{BCR98,Mar08,Sc24}; matrix versions
	were subsequently obtained by Gondard--Ribenboim \cite{GR74},
	Procesi--Schacher \cite{PS76}, and
	Djokovi\'c \cite{Djo76}. The situation is more rigid when the
	variables are freely noncommuting and positivity is tested on matrices
	of all sizes. Helton \cite{Hel02} and McCullough \cite{McC01} proved that a noncommutative
	polynomial positive semidefinite on all matrix tuples is a sum
	of hermitian squares of polynomials; see \cite{McC01} also 
	for foundational operator-valued factorization results.
	Further approaches to noncommutative
	Positivstellens\"atze, based on separation arguments,
	quadratic modules and operator algebras, were
	developed in \cite{Oza13,NT13,ANT19}.
	
	The convex counterpart concerns positivity on matricial solution sets
	of linear pencils, the so-called free spectrahedra \cite{EPS24}. 
	For linear \emph{matrix} pencils, 
	\cite{HKM} gives a convex Positivstellensatz with optimal degree bounds.
	Linear \emph{operator} inequalities provide the natural framework for working with
	arbitrary closed matrix convex sets, under the usual interior
	normalization, by the matricial Hahn--Banach theorem of
	Effros--Winkler \cite{EW} (cf.~\cite{HM12}). 
	The interplay between matrix convexity, operator
	systems, spectrahedral containment and dilation
	theory has been studied extensively; see, e.g.,
	\cite{FNT17,PSS18} and the recent survey \cite{EPS24}.
	Operator-pencil Positivstellens\"atze
	were studied by Zalar \cite{Zal17} (see also \cite{DDOSS}); more recently,
	\cite{JKM25} established operator-valued polynomial
	factorization and an operator-valued convex
	Positivstellensatz \cite{JKMconvex}. These results motivate the passage from polynomial to rational functions with  operator coefficients.
	
	The rational setting is more delicate because evaluations may be
	undefined, and the singularities are encoded by linear realizations
	of rational expressions. Klep, Pascoe and Vol\v{c}i\v{c} \cite{KPV}
	established a sum-of-squares theorem for positive
	noncommutative rational functions without hermitian
	singularities. Subsequently, Vol\v{c}i\v{c} \cite{Vol}
	treated arbitrary rational singularities and
	developed a rational convex Positivstellensatz
	using extensions of multiplication operators;
	see also \cite{Pas18}. We establish a corresponding theorem with 
	operator coefficients and an explicit complexity bound independent of
	the coefficient Hilbert space. An additional ingredient is a corrected
	domain-preserving extension mechanism; we explain its relation to
	\cite{Vol} in Subsection \ref{ssec:OG} below.
	
	\subsection{Notation and setting}
	Noncommutative (nc) rational functions arise by
	applying addition, multiplication and inversion to
	elements of the free algebra
	$\C\langle x_1,\ldots,x_{\gen}\rangle$.
	More precisely, they are equivalence classes of
	rational expressions under agreement of their
	matrix evaluations, and form the free skew field,
	the universal skew field of fractions of the free
	algebra. Its construction goes back to Amitsur
	\cite{Ami66}, while its algebraic structure was
	developed systematically by Cohn \cite{Coh}.
	The viewpoint of matrix evaluations and realizations,
	which connects nc rational functions to free analysis,
	systems theory and free real algebraic geometry,
	was developed in \cite{HMV,KVV2}; see also
	\cite{Vol2}.
	
	A complementary analytic approach to noncommutative rational
	functions arises from their action as multipliers
	on the full Fock space. Jury, Martin and Shamovich
	studied noncommutative rational multipliers and their
	factorizations \cite{JMS21a,JMS21b}, as well as
	the associated rational Clark measures
	\cite{JMS23}. More recently, Augat, Martin and
	Shamovich developed operator realization methods
	for more general noncommutative analytic functions, including
	classes of meromorphic functions \cite{AMS25}.

	We work over $\C$ with freely noncommuting, self-adjoint variables
	$x=(x_1,\ldots,x_\gen)$. Write $\rx$ for the free skew field 
	of noncommutative rational functions
	and
	$\rex$ for the nondegenerate formal rational expressions. An expression
	$r\in\rex$ determines a rational function $\rr\in\rx$, and its
	\emph{hermitian domain} $\hdom r$ consists of the hermitian matrix
	tuples at which all the inverses occurring in $r$ can be evaluated.
	We consider finite algebraic sums
	\[
	\mF=\sum_{a=1}^N C_a\otimes\rr_a
	\in\cB(\cK)\otimes\rx,
	\]
	where $\cK$ is a nonzero complex Hilbert space and $\mF=\mF^*$.
	Throughout, positivity of $\mF$ refers to all matrix sizes for
	which the indicated rational expressions are defined.
	
	Let $\cH$ be a nonzero complex Hilbert space, and fix a monic
	hermitian linear operator pencil
	\[
	L=I_{\cH}+\sum_{j=1}^{\gen}H_jx_j,
	\qquad H_j=H_j^*\in\cB(\cH).
	\]
	Its positivity set and its coefficient operator system are
	\[
	\cD(L)=\bigcup_{n\ge1}
	\bigl\{X\in\her{n}^{\gen}:L(X)\succeq0\bigr\},
	\qquad
	\cS_L=\spa_\C\{I_{\cH},H_1,\ldots,H_\gen\}.
	\]
	The set $\cD(L)$ is matrix convex. Allowing an infinite-dimensional
	$\cH$ is essential here: by \cite{HM12,EW}, linear operator inequalities
	also describe general closed matrix convex sets after the usual
	interior-point normalization.
	
	For a finite admissible set $R\subset\rex$, let $\cR$ be the
	corresponding set of rational functions and put
	\[
	V_0=\C1,\qquad
	V_\ell=\spa_\C\{\rr_1\cdots\rr_\ell:
	\rr_i\in\cR\}\quad(\ell\ge1).
	\]
	Here admissibility means that $R$ contains $1,x_1,\ldots,x_\gen$,
	is closed under adjoints and subexpressions, and contains a rational
	expression having as denominator each irreducible component of a
	denominator of an expression in $R$; see Section~\ref{s:mult} for
	details. Every finite family of rational expressions is contained in
	such a set. We measure the complexity of expressions by
	\[
	\begin{aligned}
		\tau(\alpha)&=0\quad(\alpha\in\C),&
		\tau(x_j)&=1,\\
		\tau(r+s)&=\max\{\tau(r),\tau(s)\},&
		\tau(rs)&=\tau(r)+\tau(s),\\
		\tau(r^{-1})&=2\tau(r),&
		\tau(R)&=\max_{r\in R}\tau(r).
	\end{aligned}
	\]
	The dependence on expressions, rather than only on the rational
	functions they represent, is intentional.
	
	\subsection{Main results}
	
	Our main result is a convex rational Positivstellensatz
	for matrix- and operator-valued rational functions.
	It has the same complexity bound in the
	matrix-valued and infinite-dimensional cases.
	
	\begin{thmA}\label{t:intro-main}
		Let $R\subset\rex$ be finite and admissible, and let
		$\mF=\sum_{a=1}^N C_a\otimes\rr_a
		\in\cB(\cK)\otimes V_1$
		be hermitian, where $r_a\in R$ represents $\rr_a$.
		Let $L=I_{\cH}+\sum_jH_jx_j$ be a monic hermitian
		linear operator pencil, and set
		$\ell_*=2\tau(R)+1$.
		
		Then
		\begin{equation}\label{e:intro-positive}
			\mF(X)\succeq0
			\qquad
			\text{for every }
			X\in\cD(L)\cap\bigcap_{r\in R}\hdom r
		\end{equation}
		if and only if the following certificate exists:
		\begin{enumerate}[\rm (a)]
			\item\label{thm:it:a} If $\cK$ is infinite-dimensional, there are
			$\rr,\qq\in\cB(\cK)\otimes V_{\ell_*}$ and a unital
			completely positive map $\pi:\cS_L\to\cB(\cK)$ such that
			\begin{equation}\label{e:intro-opcert}
				\mF=\rr^*\rr+\qq^*\pi(L)\qq.
			\end{equation}
			\item\label{thm:it:b} If $\cK=\C^d$, writing $m=\dim V_{\ell_*}$,
			there are $\rr,\qq\in\opm_{md\times d}(V_{\ell_*})$
			and a unital completely positive map
			$\pi:\cS_L\to\opm_{md}(\C)$ satisfying
			\begin{equation}\label{e:intro-matcert}
				\mF=\rr^*\rr+\qq^*\pi(L)\qq.
			\end{equation}
		\end{enumerate}
		Here $\pi(L)=I+\sum_j\pi(H_j)x_j$ in the respective
		coefficient algebra. In particular, the truncation bound $\ell_*$
		is independent of $\dim\cK$.
	\end{thmA}

	Part~\eqref{thm:it:b} is Theorem~\ref{t:finite-fixed}, expressed through
	Proposition~\ref{p:finite-gram}, while part~\eqref{thm:it:a} is
	Theorem~\ref{t:operator-K}. Notice that neither the size of the
	matrices on which positivity is tested nor boundedness of $\cD(L)$
	enters the hypothesis. The same truncation bound is valid for all
	finite-dimensional compressions of the coefficient operators; this
	uniformity is crucial for the operator-valued conclusion.
	
	When $\cH$ is finite-dimensional, the certificate
	in Theorem~\ref{t:intro-main} can be expressed
	directly in terms of $L$, without the unital completely
	positive (UCP) map. This gives the following two corollaries.
	
	\begin{cor}	
		Suppose that $\cH=\C^h$ and $\cK=\C^d$, and let
		$R$, $L$ and $\mF$ be as in Theorem~\ref{t:intro-main}.
		Write $m=\dim V_{\ell_*}$.
		Then \eqref{e:intro-positive} holds if and only if
		there exist
		$\rr\in\opm_{md\times d}(V_{\ell_*})$ and
		$\qq_1,\ldots,\qq_s\in
		\opm_{h\times d}(V_{\ell_*})$, with $s\le hmd$,
		such that
		\[
		\mF=\rr^*\rr+\sum_{\nu=1}^s\qq_\nu^*L\qq_\nu.
		\]
	\end{cor}
	
	\begin{proof}
		By Theorem~\ref{t:intro-main}, we have
		$\mF=\rr^*\rr+\qq^*\pi(L)\qq$,
		where $\pi:\cS_L\to\opm_{md}(\C)$ is UCP.
		By Arveson's extension theorem, $\pi$ extends
		to a UCP map on $\opm_h(\C)$, which admits
		a Kraus representation
		\[
		\pi(B)=\sum_{\nu=1}^s V_\nu^*BV_\nu,
		\qquad
		V_\nu\in\opm_{h\times md}(\C),
		\quad s\le hmd.
		\]
		Setting $\qq_\nu=V_\nu\qq$ gives the
		required representation. The converse follows
		from the positivity of $L(X)$ on $\cD(L)$.
	\end{proof}
	
	\begin{cor}
		Suppose that $\cH$ is finite-dimensional and
		$\cK$ is infinite-dimensional, and let $R$, $L$
		and $\mF$ be as in Theorem~\ref{t:intro-main}.
		Let $\cE$ be a Hilbert space such that
		$\cK$ is unitarily isomorphic to $\cE\otimes\cH$.
		Then \eqref{e:intro-positive} holds if and only if
		there exist
		$ \rr,\qq\in
		\cB(\cK,\cE\otimes\cH)\otimes V_{\ell_*} $
		such that
		\[
		\mF=\rr^*\rr+\qq^*(I_{\cE}\otimes L)\qq.
		\]
	\end{cor}
	
	\begin{proof}
		By Theorem~\ref{t:intro-main}, we have
		\[
		\mF=\rr_0^*\rr_0+\qq_0^*\pi(L)\qq_0,
		\]
		where $\rr_0,\qq_0\in\cB(\cK)\otimes V_{\ell_*}$
		and $\pi:\cS_L\to\cB(\cK)$ is UCP.
		
		Since $\cH$ is finite-dimensional, Arveson's
		extension theorem and Stinespring's theorem give
		an isometry $V:\cK\to\cE\otimes\cH$ such that
		\[
		\pi(B)=V^*(I_{\cE}\otimes B)V,
		\qquad B\in\cS_L.
		\]
		Here we use the infinite-dimensionality of $\cK$
		to choose $\cE$ large enough for the Stinespring
		representation.
		
		Choose a unitary $U:\cK\to\cE\otimes\cH$
		and set $\rr=U\rr_0$ and $\qq=V\qq_0$.
		Then
		\[
		\mF=\rr^*\rr+\qq^*(I_{\cE}\otimes L)\qq,
		\]
		as required. The converse follows from the
		positivity of $L(X)$ on $\cD(L)$.
	\end{proof}

	Taking $L=I$ recovers the noncommutative rational Hilbert's 17th
	problem in its operator-valued form: if $\mF$ is positive
	semidefinite wherever the chosen rational expressions are defined,
	then it is a finite sum of hermitian rational squares. In fact, if
	$\dim\cK=\infty$, then one square suffices,
	\begin{equation}\label{e:intro-H17}
		\mF=\rs^*\rs,
		\qquad \rs\in\cB(\cK)\otimes\rx.
	\end{equation}
	If $\dim\cK<\infty$, the number of matrix-valued squares can be
	bounded in terms of $\dim V_{\ell_*}$; see
	Corollaries~\ref{c:operator-H17} and~\ref{c:finite-H17}.
	This contrasts with the commutative operator-valued case, where
	there is no analogous rational sum-of-squares theorem in infinite
	dimensions (Section~\ref{ssec:comm-case}).

	\subsection{Relation to earlier work}\label{ssec:OG}
	The rational framework and part of the extension strategy in this
	paper originate in \cite{Vol}. In revisiting its proof, we observed
	that the hypotheses of \cite[Proposition~3.2 and Theorem~3.3]{Vol}
	do not suffice for the asserted completion conclusions. The issue
	appears for full linear pencils that are not irreducible: 
	full column rank of an
	evaluation of the pencil need not give full column rank for its
	individual irreducible components, although the latter condition is
	necessary for the proposed completion. A concrete counterexample is
	given after Proposition~\ref{p:side}.
	
	We remedy this by requiring the rank conditions for each irreducible component
	(Proposition~\ref{p:side}) and by using a more restrictive
	infinite-dimensional domain that records invertibility of all
	irreducible denominator components. These changes yield the
	finite-dimensional hermitian completion result
	(Proposition~\ref{p:main}) needed in the multiplication-operator
	argument. In particular, they supply a valid replacement for the
	extension step underlying \cite[Proposition~3.5]{Vol} and hence a
	route to the rational Positivstellensatz and Hilbert's 17th problem.
	The central approach of \cite{Vol}, that is, linear representations,
	finite-dimensional extensions, and rational multiplication
	operators, remains fundamental to the present proof. The new
	operator-valued conclusions additionally require a Gram--completely
	positive (Gram--CP) description, finite-dimensional coefficient extraction,
	and ultraweak closedness of the relevant quadratic modules.
	
	\subsection{Guide to the paper}
	Section~\ref{s:prelim} recalls free skew fields and linear
	representations. Section~\ref{s:dom} establishes the corrected
	domain-preserving extension results, and Section~\ref{s:mult}
	develops multiplication operators associated with admissible sets.
	In Section~\ref{s:cone} we introduce truncated rational quadratic
	modules and their Gram--CP representations.
	Section~\ref{s:abhay-repro} constructs reproducing tuples and
	coefficient-extraction functionals; these yield the uniform bounds
	and product-ultraweak closedness in Section~\ref{s:weighted-closed}.
	Finally, Section~\ref{s:finite-fixed} proves the matrix-valued
	rational Positivstellensatz, while Section~\ref{s:infinite-K}
	passes to general operator coefficients and derives the
	Hilbert's 17th problem corollaries, including the contrast with the
	commutative setting.

	\section{Preliminaries}\label{s:prelim}
	
	Let $\mat{m\times n}$ denote the space of complex $m\times n$ matrices, and $\mat{n}=\mat{n\times n}$. Let $\her{n}$ denote the real space of hermitian $n\times n$ matrices. For $X=(X_1,\dots,X_\gen)\in \mat{m\times n}^\gen$, $A\in\mat{p\times m}$ and $B\in\mat{n\times q}$ we write 
	$$AXB= (AX_1B,\dots,AX_\gen B)\in \mat{p\times q}^\gen, \qquad X^*=(X_1^*,\dots,X_\gen^*)\in\mat{n\times m}^\gen.$$
	
	\subsection{Free skew field}
	
	We define noncommutative rational functions using formal rational expressions and their matrix evaluations following 
	\cite{HMV} and
	\cite{KVV2}. {\bf Formal rational expressions} are syntactically valid combinations of scalars, freely noncommuting variables $x=(x_1,\dots,x_\gen)$, rational operations and parentheses. 
	A {\bf subexpression} of a formal rational expression $r$ is any formal rational expression which appears in the construction of $r$. %
	
	Given a formal rational expression $r$ and $X\in \mat{n}^\gen$, the evaluation $r(X)$ is defined in the natural way if all inverses appearing in $r$ exist at $X$. The set of all $X\in\mat{n}^\gen$ such that $r$ is defined at $X$ is denoted $\dom_n r$. The {\bf (matricial) domain} of $r$ is
	$$\dom r = \bigcup_{n\in\N} \dom_n r.$$
	Note that $\dom_n r$ is a Zariski open set in $\mat{n}^\gen$ for every $n\in\N$. A formal rational expression $r$ is {\bf nondegenerate} if $\dom r\neq\emptyset$; let $\rex$ denote the set of all nondegenerate formal rational expressions. On $\rex$ we define an equivalence relation $r_1\sim r_2$ if and only if $r_1(X)=r_2(X)$ for all $X\in\dom r_1\cap \dom r_2$. Equivalence classes with respect to this relation are called {\bf noncommutative  rational functions}. By \cite[Proposition 2.2]{KVV2} they form a skew field denoted $\rx$, which is the universal skew field of fractions of the free algebra $\px$ by \cite[Section 4.5]{Coh}. The equivalence class of $r\in\rex$ is denoted $\rr\in\rx$.
	
	There is a unique involution $*$ on $\rx$ that is determined by $\alpha^*=\overline{\alpha}$ for $\alpha\in\C$ and $x_j^*=x_j$ for $j=1,\dots,\gen$. Furthermore, this involution lifts to an involutive map $*$ on the set $\rex$.
	For $r\in\rex$ let $\hdom_n r= \dom_n r\cap \her{n}^\gen$. Then
	$$\hdom r = \bigcup_{n\in\N} \hdom_n r$$
	is the {\bf hermitian domain} of $r$. Note that $\hdom_n r$ is Zariski dense in $\dom_n r$ because $\her{n}$ is Zariski dense in $\mat{n}$ and $\dom_n r$ is Zariski open in $\mat{n}^\gen$. 
	Finally, while a noncommutative rational function $\rr\in\rx$ can be represented by various expressions, there exists an expression $r$ representing $\rr$ whose hermitian domain contains the hermitian domain of any other expression representing $\rr$ \cite[Proposition 2.1]{Vol}.

	\subsection{Linear representation of a formal rational expression}
	
	Fundamental tools for handling noncommutative rational functions are linear representations \cite{CR,HMS}. Let $r\in\rex$. By \cite[Theorem 4.2 and Algorithm 4.3]{HMS} there exist $e\in\N$, vectors $u,v\in\C^e$ and an affine matrix pencil $M=M_0+M_1x_1+\cdots+M_\gen x_\gen$ with $M_j\in\mat{e}$ satisfying the following. 
	For every unital $\C$-algebra $\cA$ and $a\in\cA^\gen$, 
	\begin{enumerate}[(i)]
		\item if $r$ can be evaluated at $a$, 
		then $M(a)=M_0\otimes I+M_1\otimes a_1+\cdots +M_\gen\otimes a_\gen \in\GL_e(\cA)$ and
		$r(a) = (u\otimes I)^* M(a)^{-1}(v\otimes I)$;
		\item if $M(a)\in\GL_e(\cA)$ and $\cA=\mat{n}$ for some $n\in\N$, then $r$ can be evaluated at $a$.
	\end{enumerate}
	We say that the triple $(u,M,v)$ is a {\bf linear representation} of $r$ of size $e$, and $M$ is a {\bf denominator} of $r$.
	An affine matrix pencil $M$ is {\bf full} \cite[Section 1.4]{Coh} if it is invertible over $\rx$, or equivalently, $\det M(X)\neq0$ for some matrix tuple $X$ (see, for example, \cite[Proposition 2.4]{Vol}). 
	In particular, full pencils are precisely those that can appear in  linear representations of nondegenerate rational expressions. 
	
	A pencil is {\bf irreducible} if it cannot be put into a $2\times2$ block upper triangular form with square diagonal blocks
	$(\begin{smallmatrix}\star & \star \\ 0 & \star\end{smallmatrix})$ 
	by a left and a right basis change.
	Every nonzero irreducible pencil is full. On the other hand, every full pencil is, up to a left and a right basis change, equal to a block upper-triangular pencil whose diagonal blocks are irreducible pencils. Given a pencil $M$, we refer to such nonconstant blocks, unique up to basis change, as {\bf irreducible components} of $M$. 
	Irreducible components of a denominator for $r$, disregarding
	multiplicities, are unique up to basis change. %
	Indeed, let $M$ and $M'$ be denominators of $r$, and suppose $M'$ is an irreducible component that is not a component of $M$ after a basis change. Then, by \cite[Theorem 2.12]{HKV} there is $X\in\mat{n}^\gen$ such that $\det M'(X)=0$ but $\det M(X)\neq0$, which contradicts the invertibility requirement in the definition of a linear representation.

	\section{Extensions and retractions}\label{s:dom}

	In this section we establish domain-preserving completion results,
	first for full linear pencils (Proposition~\ref{p:side})
	and then for hermitian tuples in rational domains
	(Proposition~\ref{p:main}). Lemma~\ref{l:change-basis-completion}
	provides a variant compatible with changes of inner product.

	For the purpose of this section we extend evaluations of linear matrix pencils to tuples of rectangular matrices.
	If $\Lambda=\sum_{j=1}^\gen \Lambda_jx_j$ is of size $e$ and $X\in\mat{\ell \times m}^\gen$ then
	$$\Lambda(X)=\sum_{j=1}^\gen \Lambda_j\otimes X_j \in \mat{e\ell\times em}.$$
	
	We require a technical tool.
	
	\begin{lem}[\protect{\cite[Lemma 3.1]{Vol}}]\label{l:irr}
		Let $\Lambda = \sum_{j=1}^\gen\Lambda_jx_j$ be a homogeneous irreducible pencil of size $e$. Let $\ell\le m$ and denote $n=(m-\ell)(e-1)$. Given $C\in\mat{me\times \ell e}$, consider the pencil $\widetilde{\Lambda}$ of size $(m+n)e$ in $\gen(m+n)(n+m-\ell)$ variables $z_{jpq}$,
		$$\widetilde{\Lambda} = 
		\begin{pmatrix}
			C & 0 \\ 0 & 0
		\end{pmatrix}
		+\sum_{j=1}^\gen\sum_{q=1}^{n+m-\ell}\left(
		\sum_{p=1}^m
		\begin{pmatrix}
			0 & \widehat{E}_{p,q}\otimes \Lambda_j \\ 0 & 0
		\end{pmatrix}z_{jpq}
		+\sum_{p=m+1}^{m+n}
		\begin{pmatrix}
			0 & 0 \\ 0 & \widecheck{E}_{p-m,q}\otimes \Lambda_j
		\end{pmatrix}z_{jpq}
		\right)
		$$
		where $\widehat{E}_{p,q} \in\mat{m\times (n+m-\ell)}$ and $\widecheck{E}_{p-m,q} \in\mat{n\times (n+m-\ell)}$ are the standard matrix units. If $C$ has full rank, then the pencil $\widetilde{\Lambda}$ is full.
	\end{lem}
	
	The following proposition strengthens the hypotheses of \cite[Proposition 3.2]{Vol}.

	\begin{prop}\label{p:side}
		Let $\Lambda$ be a homogeneous full pencil of size $e$. Suppose $X\in\mat{m\times \ell}^\gen$ with $\ell\le m$ is such that $\Lambda'(X)$ has full rank for every irreducible component $\Lambda'$ of $\Lambda$. Then there exist $\widehat{X}\in \mat{m\times (n+m-\ell)}^\gen$ and $\widecheck{X}\in \mat{n\times (n+m-\ell)}^\gen$ for some $n\in\N$ such that
		\begin{equation}\label{e:side}
			\det \Lambda
			\begin{pmatrix}
				X & \widehat{X} \\ 0 & \widecheck{X}
			\end{pmatrix}
			\neq0.
		\end{equation}
	\end{prop}
	
	\begin{proof}
		First, observe that if the conclusion holds for some $n$, it holds for all sufficiently large $n$. Indeed, for every $k\ge e-1$ there is $Z_k\in\mat{k}^\gen$ such that $\det \Lambda(Z_k)\neq0$ by \cite[Theorem 1.8]{DM}. Thus, if \eqref{e:side} holds, then also
		$$\det \Lambda
		\begin{pmatrix}
			X & \widehat{X} & 0\\ 0 & \widecheck{X}&0 \\0&0&Z_k
		\end{pmatrix}
		\neq0$$
		for every $k\geq e-1$. 
		
		Second, suppose that the proposition holds for irreducible pencils. Given a pencil $\Lambda$,  the proposition then holds for its irreducible components $\Lambda_1,\dots,\Lambda_t$.        
		By the first observation, 
		there is $n\in\N$ such that for every $i=1,\dots,t$, the set $\cO_i$ of pairs
		$(\widehat{X},\widecheck{X})\in \mat{m\times (n+m-\ell)}^\gen\times \mat{n\times (n+m-\ell)}^\gen$ satisfying \eqref{e:side} for $\Lambda_i$ is Zariski open and nonempty. 
		Hence, $\cO_1\cap\cdots\cap \cO_t\neq\emptyset$, so the proposition also holds for full pencils. Thus, we can without loss of generality assume that $\Lambda$ is irreducible.
		
		If $m=\ell$, full column rank makes $\Lambda(X)$
		invertible. A block-diagonal extension with any invertible
		evaluation of the full pencil $\Lambda$ then proves the assertion.
		Hence in the following construction we may assume $m>\ell$.
		
		Let $n_1=(m-\ell)(e-1)$ and $\widetilde e=(m+n_1)e$. Apply Lemma \ref{l:irr} to $C=\sum_{j=1}^\gen X_j\otimes \Lambda_j$; since the resulting pencil $\widetilde{\Lambda}$ is full, there exists
		$Z \in\mat{\widetilde e -1}^{\gen(m+n_1)(n_1+m-\ell)}$ such that $\widetilde{\Lambda}(Z)$ is invertible. Therefore the matrix
		$$
		\sum_{j=1}^\gen\Lambda_j\otimes\left(
		\left(\begin{smallmatrix}I\otimes X_j & 0 \\ 0 & 0\end{smallmatrix}\right)
		+\sum_{q=1}^{n_1+m-\ell}\left(\sum_{p=1}^{m}
		\left(\begin{smallmatrix}
			0 & Z_{jpq}\otimes \widehat{E}_{p,q} \\ 0 & 0
		\end{smallmatrix}\right)
		+\sum_{p=m+1}^{m+n_1}
		\left(\begin{smallmatrix}
			0 & 0 \\ 0 & Z_{jpq}\otimes \widecheck{E}_{p-m,q}
		\end{smallmatrix}\right)
		\right)\right)
		$$
		is invertible since it is similar to $\widetilde{\Lambda}(Z)$ (via a permutation matrix).
		Thus there are $\widehat{X}\in \mat{m\times (n+m-\ell)}^\gen$ and $\widecheck{X}\in \mat{n\times (n+m-\ell)}^\gen$ such that
		$$\det \Lambda
		\begin{pmatrix}
			X & \widehat{X} \\ 0 & \widecheck{X}
		\end{pmatrix}
		\neq0,$$
		where $n=(\widetilde e -2)m+n_1(\widetilde e -1)$. 
	\end{proof}
	
	\begin{rem}
		In  \cite[Proposition 3.2]{Vol}, only the full rank of $\Lambda(X)$ is imposed, instead of the full rank for all irreducible components. However, taking
		$$\Lambda=\begin{pmatrix}x_1&x_2\\0&x_3\end{pmatrix},\quad
		X=\left(
		\begin{pmatrix}1\\0\end{pmatrix},
		\begin{pmatrix}0\\1\end{pmatrix},
		\begin{pmatrix}0\\0\end{pmatrix}
		\right)
		$$
		shows that while $\Lambda(X)$ has full rank, there are no $\widehat{X},\widecheck{X}$ satisfying \eqref{e:side}.
	\end{rem}
	
	Let $\cfm$ be the algebra of $\N\times\N$ matrices over $\C$ that have only finitely many nonzero entries in each column; that is, elements of $\cfm$ can be viewed as linear operators on $\oplus^{\N}\C$. Given $r\in\rex$ with a denominator $M$, let $\dom_\infty r$ denote the set of tuples $X\in\cfm^\gen$ such that $M'(X)$ is invertible for every irreducible component $M'$ of $M$. Note that for $X\in\dom_\infty r$, both $r(X)$ and $M(X)$ are well-defined, and $r(X)=(u^*\otimes I) M(X)^{-1}(v\otimes I)$ where $(u,M,v)$ is a linear representation of $r$. While a denominator of  $r$ is not unique, the set $\dom_\infty r$ is independent of the choice of the denominator (by uniqueness of irreducible components of denominators). 
	Note that this definition of $\dom_\infty r$ is slightly more restrictive than the one in \cite[Section 3]{Vol}.
	
	\begin{prop}\label{p:main}
		Let $r_1,\dots,r_t\in\rex$. If $X\in\her{\ell}^\gen$ and $Y\in\mat{m\times \ell}^\gen$ are such that
		$$\begin{pmatrix}
			X & \begin{matrix} Y^* & 0\end{matrix} \\
			\begin{matrix} Y \\ 0\end{matrix} & \mathlarger{W}
		\end{pmatrix} 
		\in \bigcap_{i=1}^t\dom_\infty r_i
		$$
		for some $W\in\cfm^\gen$, then there exist $n\ge m$, $E\in\gl{n}$ and $Z\in\her{n}^\gen$ such that
		$$\begin{pmatrix}
			X & \begin{pmatrix} Y^* & 0\end{pmatrix}E^* \\
			E\begin{pmatrix} Y \\ 0\end{pmatrix} & \mathlarger{Z}
		\end{pmatrix} \in \bigcap_{i=1}^t\hdom r_i.
		$$
	\end{prop}
	
	\begin{proof}
		Let $M$ be the direct sum of denominators of $r_1,\dots,r_t$, of size $e$.
		By assumption,
		$$M'\begin{pmatrix}
			X & \begin{matrix} Y^* & 0\end{matrix} \\
			\begin{matrix} Y \\ 0\end{matrix} & \mathlarger{W}
		\end{pmatrix} 
		$$
		is an invertible matrix over $\cfm$ for every irreducible component $M'$ of $M$. If $M'=M_0'+M_1'x_1+\cdots+M_\gen'x_\gen$, then the matrices
		$$M_0'\otimes \begin{pmatrix}I \\ 0\end{pmatrix}+
		\sum_{j=1}^\gen M_j'\otimes \begin{pmatrix}X_j \\ Y_j\end{pmatrix},\qquad
		M_0'\otimes \begin{pmatrix}I & 0\end{pmatrix}+
		\sum_{j=1}^\gen M_j'\otimes \begin{pmatrix}X_j & Y_j^*\end{pmatrix}
		$$
		have full rank. 
		
		Let $\Lambda$ be the homogenization of $M$ with $x_0$. By Proposition \ref{p:side} and its transpose there exist $k\in\N$ and
		\begin{align*}
			A'\in \mat{\ell\times (k+m)}^{1+\gen},\quad 
			B'\in \mat{m\times (k+m)}^{1+\gen},\quad
			C'\in \mat{k\times (k+m)}^{1+\gen}, \\
			A''\in \mat{(k+m)\times \ell}^{1+\gen},\quad 
			B''\in \mat{(k+m)\times m}^{1+\gen},\quad
			C''\in \mat{(k+m)\times k}^{1+\gen} 
		\end{align*}
		such that
		$$\Lambda
		\begin{pmatrix}
			X & A' \\ Y & B' \\ 0 & C'
		\end{pmatrix}
		,\qquad
		\Lambda
		\begin{pmatrix}
			X & Y^* &0 \\ A'' & B'' & C''
		\end{pmatrix}
		$$
		are invertible; here, $X_0=I$ and $Y_0=0$. 
		Consequently, there exists $\ve\in\C\setminus\{0\}$ such that
		$$\begin{pmatrix}
			\Lambda(X) & 0 & 0 &0 & 0 \\
			0 & 0 & 0 &0 & 0 \\
			0 & 0 & 0 &0 & 0 \\
			0 & 0 & 0 &0 & 0 \\
			0 & 0 & 0 &0 & 0
		\end{pmatrix}
		+\ve\begin{pmatrix}
			\begin{pmatrix}0 & 0 \\ 0 &0 \end{pmatrix} &
			\Lambda
			\begin{pmatrix}
				X & Y^* &0 \\ A'' & B'' & C''
			\end{pmatrix}
			\\
			\Lambda
			\begin{pmatrix}
				X & A' \\ Y & B' \\ 0 & C'
			\end{pmatrix}
			&
			\begin{pmatrix}0 & 0 & 0\\ 0 &0 &0\\ 0 &0 &0\end{pmatrix}
		\end{pmatrix}$$
		is invertible; this matrix is similar to
		\begin{equation}\label{e:eps}
			\Lambda
			\begin{pmatrix}
				X & 0 & \ve X & \ve Y^* &0 \\
				0 & 0 & \ve A''& \ve B'' & \ve C'' \\
				\ve X  & \ve A' & 0 &0 &0 \\
				\ve Y & \ve B' & 0 &0 &0 \\
				0 & \ve C' & 0 &0 &0
			\end{pmatrix}.
		\end{equation}
		Thus the matrix \eqref{e:eps} is invertible; its block structure and the linearity of $\Lambda$ imply that \eqref{e:eps} is invertible for every $\ve\neq0$, so we can choose $\ve=1$. After performing elementary row and column operations on \eqref{e:eps} we conclude that
		$$\Lambda
		\begin{pmatrix}
			X & Y^* & 0 & 0 &0 \\
			Y & 0 & - Y & B' &0 \\
			0  & - Y^* & - X &  A' &0 \\
			0 & B'' &  A''& 0 & C'' \\
			0 & 0 & 0 & C' &0
		\end{pmatrix}
		$$
		is invertible.
		
		We have seen that for $n=2(m+k)+\ell$ there exists $Z'\in \mat{n}^{1+\gen}$ such that
		\begin{equation}\label{e:2psd}
			\det\left(M_0\otimes
			\begin{pmatrix}
				I & \begin{matrix} 0 & 0\end{matrix} \\
				\begin{matrix} 0 \\ 0\end{matrix} & \mathlarger{Z'_0}
			\end{pmatrix}+
			\sum_{j=1}^\gen M_j\otimes 
			\begin{pmatrix}
				X_j & \begin{matrix} Y_j^* & 0\end{matrix} \\
				\begin{matrix} Y_j \\ 0\end{matrix} & \mathlarger{Z'_j}
			\end{pmatrix}
			\right)\neq0.
		\end{equation}
		The set of all $Z'\in\mat{n}^{1+\gen}$ satisfying \eqref{e:2psd} is thus a nonempty Zariski open set in $\mat{n}^{1+\gen}$. Since the set of positive definite $n\times n$ matrices is Zariski dense in $\mat{n}$, there exists $Z'\in\her{n}^{1+\gen}$ with $Z'_0\succ0$ such that \eqref{e:2psd} holds. If $Z_0' = E^{-1}E^{-*}$, let $Z_j=EZ'_jE^*$ for $1\le j \le \gen$. Then
		$$M\begin{pmatrix}
			X & \begin{pmatrix} Y^* & 0\end{pmatrix}E^* \\
			E\begin{pmatrix} Y \\ 0\end{pmatrix} & \mathlarger{Z}
		\end{pmatrix}
		$$
		is invertible, so
		$$
		\begin{pmatrix}
			X & \begin{pmatrix} Y^* & 0\end{pmatrix}E^* \\
			E\begin{pmatrix} Y \\ 0\end{pmatrix} & \mathlarger{Z}
		\end{pmatrix}  \in \bigcap_{i=1}^t\hdom r_i
		$$
		by the definition of linear representations.
	\end{proof}
	
	The following lemma extends the result of Proposition~\ref{p:main} to tuples that are not necessarily hermitian but are simultaneously similar to a hermitian tuple in the prescribed rational domains.

	\begin{lem}	\label{l:change-basis-completion}
		Let $r_1,\ldots,r_t\in\rex$, and let $S,U_0,W_0$ be
		finite-dimensional Hilbert spaces. Suppose that
		$T=(T_1,\ldots,T_\gen)$ is a tuple of operators on
		$S\oplus U_0\oplus W_0$ of the form
		\[
		T_j=
		\begin{pmatrix}
			X_j & * & *\\
			Y_j & * & *\\
			0   & * & *
		\end{pmatrix},
		\qquad X_j=X_j^*,\quad 1\le j\le\gen.
		\]
		Assume that $T$ is simultaneously similar to a hermitian
		tuple $\widehat X\in\bigcap_{i=1}^t\hdom r_i$.
		Then there exist $k\in\N$, an invertible operator
		$E$ on $U_0\oplus\C^k$, and a tuple of hermitian
		operators $Z=(Z_1,\ldots,Z_\gen)$ on $U_0\oplus\C^k$
		such that
		\[
		\widetilde X_j=
		\begin{pmatrix}
			X_j &
			\begin{pmatrix}Y_j^*&0\end{pmatrix}E^*\\[2mm]
			E\begin{pmatrix}Y_j\\0\end{pmatrix} & Z_j
		\end{pmatrix},
		\qquad 1\le j\le\gen,
		\]
		satisfies
		\[
		\widetilde X\in\bigcap_{i=1}^t\hdom r_i.
		\]
	\end{lem}
	
	\begin{proof}
		Let $M'=M'_0+\sum_{j=1}^\gen M'_jx_j$ be an
		irreducible component of a denominator of one of
		the expressions $r_1,\ldots,r_t$ of size $e$. 
		Since $T$ is simultaneously similar to $\widehat X$, the
		matrix $M'(T)$ is invertible. Thus, its
		restriction to $\C^e\otimes S$ has full column rank. This restriction
		is represented by
		\[
		M'_0\otimes
		\begin{pmatrix}I_S\\0\\0\end{pmatrix}
		+
		\sum_{j=1}^\gen M'_j\otimes
		\begin{pmatrix}X_j\\Y_j\\0\end{pmatrix}.
		\]
		Moreover, $T^*$ is also simultaneously similar to
		$\widehat X$. Indeed, if $T_j = F^{-1} \widehat X_j F$, then
		$ T_j^*=F^* \widehat X_j (F^*)^{-1}. $
		Thus $M'(T^*)$ is invertible as well. Since
		$X_j=X_j^*$, its first block row is
		\[
		M'_0\otimes
		\begin{pmatrix}I_S&0&0\end{pmatrix}
		+
		\sum_{j=1}^\gen M'_j\otimes
		\begin{pmatrix}X_j&Y_j^*&0\end{pmatrix},
		\]
		which therefore has full row rank. 		
		
		The rest of the proof is the same as in the proof of Proposition
		\ref{p:main}.
	\end{proof}
	
	\section{Multiplication operators of admissible sets}\label{s:mult}

	In this section we introduce admissible sets of rational expressions
	and the associated multiplication operators.
	Lemma~\ref{l:fX} shows that these operators lie in
	the prescribed rational domains, while
	Proposition~\ref{p:gns} relates their finite-dimensional
	extensions to evaluations of rational expressions.

	We say that a finite set $R\subset\rex$ is {\bf admissible} if
	\begin{enumerate}
		\item $R$ contains $1,x_1,\dots,x_\gen$;
		\item $R$ is closed under subexpressions and $*$;
		\item if $r\in R$ has a denominator $M$ with an irreducible component $M'$, then there is $r'\in R$ with $M'$ as a denominator. 
	\end{enumerate} 
	Note that every $r\in\rex$ belongs to a finite admissible set because only finitely many (up to basis change) irreducible pencils are components of denominators of $r$. 
	Indeed, by the inclusion of domains and \cite[Theorem 2.12]{HKV}, taking subexpressions does not introduce new irreducible denominator components, and taking adjoints only produces their adjoints.
	
	Let $\cR\subset\rx$ be the set of noncommutative rational functions represented by an admissible $R$. 
	For $\ell\in\N$ we define an increasing sequence of finite-dimensional vector subspaces 
	$$V_\ell = \spa_\C\overbrace{\cR\cdots \cR}^\ell\subset \rx.$$
	Then, $V:=\bigcup_{\ell\in\N} V_\ell$ is a finitely generated $*$-subalgebra of $\rx$. For $j=1,\dots,\gen$ we define left multiplication operators
	$$\fX_j:V\to V,\qquad \fX_j \rr = x_j\rr.$$
	By \cite[Lemma 4.1]{Vol}, there is a linear functional $\phi:V\to\C$ such that $\phi(\rr^*)=\overline{\phi(\rr)}$ and $\phi(\rr\rr^*)>0$ for all $\rr\in V\setminus\{0\}$.
	This functional gives rise to an inner product on $V$,
	\begin{equation}\label{e:inner}
		(\rr_1,\rr_2)=\phi(\rr_2^*\rr_1).
	\end{equation}
	With respect to this inner product we can inductively build an ordered orthogonal basis $\mathscr{B}$ of $V$ with the property that $\mathscr{B}\cap V_\ell$ is a basis of $V_\ell$ for every $\ell\in\N$. 
	The following lemma gives rise to tuples $\fX\in\dom_\infty r$ to which Proposition \ref{p:main} will later be applied.
	
	\begin{lem}\label{l:fX}
		With respect to the inner product \eqref{e:inner} and the ordered basis $\mathscr{B}$ as above, the operators $\fX_1,\dots,\fX_\gen$ are represented by hermitian matrices in $\cfm$, and $\fX\in\dom_\infty r$ for every $r\in R$.
	\end{lem}

	\begin{proof}
		Since
		$$\big(\fX_j\rr_1,\rr_2\big)
		= \phi\big(\rr_2^*x_j\rr_1\big)
		= \big( \rr_1,\fX_j\rr_2\big)$$
		for all $\rr_1,\rr_2\in V$ 
		and $\fX_j(V_\ell)\subseteq V_{\ell+1}$ for all $\ell\in\N$, it follows that the matrix representation of $\fX_j$ with respect to $\mathscr{B}$ is hermitian and has only finitely many nonzero entries in each column and row.
		The rest follows by induction on the construction of $r$ since 
		$\fX_j$ are the left multiplication operators on $V$ and $R$ is admissible.
	\end{proof}

	We use a complexity measure $\tau:\rex\to\N\cup\{0\}$ defined recursively as follows:
	\begin{enumerate}[(i)]
		\item $\tau(\alpha)=0$ for $\alpha\in\C$ and $\tau(x_j)=1$ for $1\le j\le \gen$;
		\item $\tau(r_1+r_2)=\max\{\tau(r_1),\tau(r_2)\}$ for $r_1,r_2\in\rex$;
		\item $\tau(r_1r_2)=\tau(r_1)+\tau(r_2)$ for $r_1,r_2\in\rex$;
		\item $\tau(r^{-1})=2\tau(r)$ for $r,r^{-1}\in\rex$.
	\end{enumerate}
	For a finite set $R\subset\rex$ let $\tau(R)=\max_{r\in R}\tau(r)$.
	If $R$ is admissible, then $x_1\in R$, so $\tau(R)\ge1$.
	
	\begin{prop}[\protect{\cite[Proposition 4.3]{Vol}}]\label{p:gns}
		Let the notation be as above, and $\ell\in\N$.
		Let $r\in R$ and $s\in R\cdots R$ ($\ell$ times) satisfy $2\tau(r)+\tau(s)\le \ell+2$.	
		If $U$ is a finite-dimensional Hilbert space containing $V_{\ell+1}$ and $X$ is a $\gen$-tuple of hermitian operators on $U$ such that $X\in\hdom r$ and
		$$X_j|_{V_\ell}=\fX_j|_{V_\ell}\quad \text{for }j=1,\dots,\gen,$$
		then
		$$r(X)\rs=\rr\rs.$$
	\end{prop}

	The proposition also holds with $V_k$ replaced by $V_k^{1\times d}$, with $\fX_j$ acting componentwise and $s$ replaced by a standard unit row. Indeed, the proof of \cite[Proposition 4.3]{Vol} applies componentwise and uses only the stated multiplication identities.
	
	\section{Truncated \ratQM }
	\label{s:cone}

	We introduce truncated rational quadratic modules
	with matrix coefficients. Their Gram--CP
	description in Proposition~\ref{p:finite-gram} establishes
	convexity and provides the foundation for the later
	operator-valued extension.

	Let $\cH$ be a nonzero complex Hilbert space and let
	\[
	L=I_{\cH}+\sum_{j=1}^{\gen}H_jx_j,\qquad
	H_j=H_j^*\in\cB(\cH).
	\]
	Set $\cS_L=\operatorname{span}_{\C}\{I_{\cH},H_1,\ldots,H_\gen\}$
	and $\cD_n(L)=\{X\in\her{n}^{\gen}:L(X)\succeq0\}$.
	We write $\cD(L)=\bigcup_{n\ge1}\cD_n(L)$.
	Fix a finite admissible set $R\subset\rex$
	and let $\cR\subset\rx$ be the corresponding set of noncommutative rational functions. Let $V_k=V_k(\cR)$ be the vector spaces from Section~\ref{s:mult}. 
	For $\ell\ge0$ put
	$m=m_\ell=\dim_{\C}V_\ell$ (with $V_0=\C1$).
	
	\begin{defn}
		\label{d:matrix-corrected}
		For a positive integer $d$ define 
		\begin{equation}\label{e:matrix-corrected}
			Q_\ell^{(d)}(L;\cR )=\{\rr^*\rr+\qq^*\pi(L)\qq \ :\  
			\rr,\qq\in\opm_{m d \times d}(V_\ell), \ 
			\pi\in\UCP(\cS_L,\opm_{md}(\C))\},
		\end{equation}
		where
		$\pi(L)=I_{md}+\sum_j\pi(H_j)x_j$.
		This is a subset of 
		\[\opm_d(V_{2\ell+1})_{\rm sa}
		=\{\mF \in\opm_d(V_{2\ell+1}) : \mF=\mF^*\}.
		\]
	\end{defn}
	
	\begin{prop}\label{p:finite-gram}
		Let $\ee_1,\ldots,\ee_m$ be a basis of $V_\ell$, and set
		$ \mathbf e_\ell=\operatorname{col}(\ee_1 I_d,\ldots,\ee_m I_d). $
		Then $\mF\in \opm_d(V_{2\ell+1})_{\rm sa}$ is in 
		$Q_\ell^{(d)}(L;\cR )$ if and only if there exist
		positive semidefinite
		$G\in\opm_m(\opm_d(\C))_+$ and a completely positive map
		$\psi:\cS_L\to\opm_m(\opm_d(\C))$ satisfying
		\begin{equation}\label{e:finite-gram}
			\mF=\mathbf e_\ell^*\bigl(G+\psi(L)\bigr)\mathbf e_\ell.
		\end{equation}
		In particular, $Q_\ell^{(d)}(L;\cR )$ is a convex cone.
	\end{prop}
	\begin{proof}
		If $\mF = \rr^*\rr + \qq^*\pi(L) \qq$, write $\rr= \mathsf R\,\mathbf e_\ell$ and $\qq=\mathsf Q\,\mathbf e_\ell$ for
		$\mathsf R,\mathsf Q\in\opm_{m}( \opm_{d}(\C))$. 
		Then $G=\mathsf R^*\mathsf R$ and
		$\psi(B) = \mathsf Q^*\pi(B)\mathsf Q$ give \eqref{e:finite-gram}.

		Conversely, suppose that \eqref{e:finite-gram} holds.
		Let $\mathsf R=G^{1/2}$ and $D=\psi(I_{\cH})^{1/2}$.
		By Arveson's extension theorem and Stinespring's theorem \cite{Pau},
		there is a representation $\rho$ of $\cB(\cH)$ and an
		operator $\Gamma$ with $\psi(B) = \Gamma^*\rho(B) \Gamma$.
		Writing $\Gamma = \Theta D$ in the polar decomposition, we have
		$\Theta^* \Theta = P$, where $P$ is the support projection of $D$.
		Fix a state $\omega$ on $\cS_L$. Then
		\[
		\pi(B)=\Theta^*\rho(B) \Theta+\omega(B)(I_{md}-P)
		\]
		is a unital completely positive map and satisfies
		$\psi(B)=D\pi(B)D$ for $B\in\cS_L$.
		Taking $\rr = \mathsf R \mathbf e_\ell$ and $\qq=D\mathbf e_\ell$
		gives $\mF = \rr^*\rr + \qq^*\pi(L)\qq$, as required.
	\end{proof}

	\section{Finite-dimensional reproducing tuples and coefficient extraction}\label{s:abhay-repro}

	In this section we construct finite-dimensional hermitian tuples
	that reproduce prescribed rational functions
	(Proposition~\ref{p:reproducing-tuple}), using the
	domain-preserving extensions of
	Proposition~\ref{p:finite-extension}.
	These tuples yield coefficient-extraction functionals
	(Lemma~\ref{lem:coefficient-extraction}), which will
	be used to control the truncated quadratic modules.

	We continue with the 
	notation of Section \ref{s:cone}.

	\begin{prop}\label{p:finite-extension}
		For every $\ell\in\N$, there exist an integer
		$n_\ell\ge0$ and a $\gen$-tuple of hermitian
		operators $\widetilde X^{(\ell)}$ on
		$V_{\ell+1}\oplus\C^{n_\ell}$ such that
		\[
		\widetilde X^{(\ell)}
		\in\bigcap_{r\in R}\hdom r
		\qquad\text{and}\qquad
		\widetilde X_j^{(\ell)}|_{V_\ell}
		=\fX_j|_{V_\ell},
		\quad j=1,\ldots,\gen.
		\]
	\end{prop}
	
	\begin{proof}
		Fix $\ell\in\N$, and let $U_0$ be the orthogonal complement
		of $V_{\ell+1}$ in $V_{\ell+2}$. With respect to the
		orthogonal decomposition $V=V_{\ell+1}\oplus U_0\oplus
		V_{\ell+2}^{\perp}$, the self-adjoint multiplication tuple
		has the block form
		\[
		\fX=\begin{pmatrix}
			X & \begin{pmatrix}Y^*&0\end{pmatrix}\\
			\begin{pmatrix}Y\\0\end{pmatrix} & W
		\end{pmatrix}.
		\]
		By Lemma~\ref{l:fX}, the tuple $\fX$ belongs to
		$\bigcap_{r\in R}\dom_\infty r$.
		Applying
		Proposition~\ref{p:main} to the finitely many expressions in
		$R$, we obtain $E,Z$ such that
		\[
		\widetilde X^{(\ell)}=
		\begin{pmatrix}
			X & \begin{pmatrix}Y^*&0\end{pmatrix}E^*\\
			E\begin{pmatrix}Y\\0\end{pmatrix}&Z
		\end{pmatrix}
		\in\bigcap_{r\in R}\hdom r.
		\]
		Since $\fX_j(V_\ell)\subseteq V_{\ell+1}$,
		$\widetilde X_j^{(\ell)}|_{V_\ell}=\fX_j|_{V_\ell}$.
	\end{proof}

	\begin{prop}\label{p:reproducing-tuple}
		For every $\ell\in\N$, there exist integers
		$K_\ell\ge\ell$ and $n_{K_{\ell}} \geq 0$, and a
		$\gen$-tuple of hermitian operators $\mathfrak A^{(\ell)}$ acting on
		$ U^{(\ell)}=V_{K_\ell+1}\oplus\C^{n_{K_\ell}}, $
		such that
		$ \mathfrak A^{(\ell)}\in\bigcap_{r\in R}\hdom r $
		and
		\begin{equation}\label{eq:A-reproduces}
			\qq(\mathfrak A^{(\ell)})\one=\qq,
			\qquad \qq\in V_\ell.
		\end{equation}
	\end{prop}
	
	\begin{proof}
		Fix $\ell \in \mathbb N.$ 
		Choose an integer $K_\ell$ such that
		$ K_\ell\ge \ell $
		and
		$ K_\ell+2\ge (\ell+1) \tau(R). $ By Proposition~\ref{p:finite-extension}, the tuple
		$ \mathfrak A^{(\ell)}:=\widetilde X^{(K_\ell)} $  
		acts on $V_{K_\ell+1}\oplus\C^{n_{K_\ell}}$ and satisfies
		$\mathfrak A^{(\ell)}\in\bigcap_{r\in R}\hdom r$.
		
		\smallskip
		
		Now we shall prove \eqref{eq:A-reproduces}. Since $V_\ell$ is spanned by
		products of at most $\ell$ elements of $\cR$, it suffices by linearity to
		consider
		\[
		\qq=\qq_1\cdots\qq_k,
		\qquad 1\le k\le\ell,
		\]
		where each $\qq_j\in\cR$. Choose representatives $q_j\in R$ of
		$\qq_j$. For $1\le j\le k$, put
		$ s_j:=q_{j+1}\cdots q_k, $
		with $s_k=1$. Then
		$ \tau(s_j)\le (k-j) \tau(R), $
		and therefore
		\[
		2\tau(q_j)+\tau(s_j)
		\le 2 \tau(R) + (k-j) \tau(R)
		\le (\ell+1) \tau(R)
		\le K_\ell+2.
		\]
		Because $K_\ell\ge\ell$, each $s_j$ may be viewed as a product of
		$K_\ell$ elements of $R$ by inserting copies of $1$. Therefore, a successive application of Proposition
		\ref{p:gns} gives
		\[
		(q_1\cdots q_k)(\mathfrak A^{(\ell)})\one
		=
		\qq_1\cdots\qq_k.
		\]
		Thus \eqref{eq:A-reproduces} holds on a spanning set of $V_\ell$, and hence
		on all of $V_\ell$.
	\end{proof}
	
	\begin{lem}\label{lem:coefficient-extraction}
		Let $k\in\N$, $\cW\subseteq V_k$ be a finite-dimensional subspace, and let
		$\{\ee_1,\dots,\ee_w\}$ be an orthonormal basis of $\cW$. For each $1\le j\le w$
		define a linear functional 
		\[
		\Omega_j:\cB(U^{(k)})\longrightarrow \C,
		\qquad
		T \longmapsto \langle T\one,\ee_j \rangle.
		\]
		Let $\cK$ be a complex Hilbert space. If
		$ \mF =\sum_{j=1}^w C_j\otimes\ee_j\in\cB(\cK)\otimes \cW, $
		then
		\begin{equation}\label{eq:general-extraction}
			(\id_{\cB(\cK)}\otimes\Omega_j) \mF (\mathfrak A^{(k)}) = C_j,
			\qquad 1\le j\le w,
		\end{equation}
		where the tuple $\mathfrak A^{(k)}$ is from Proposition \ref{p:reproducing-tuple}.
	\end{lem}
	
	\begin{proof}
		By Proposition \ref{p:reproducing-tuple},
		$\ee_i( \mathfrak A^{(k)})\one=\ee_i$ for every $i$. Hence
		\[
		(\id_{\cB(\cK)}\otimes\Omega_j) \mF( \mathfrak A^{(k)})
		=
		\sum_{i=1}^w C_i\,\Omega_j(\ee_i(\mathfrak A^{(k)}))
		=
		\sum_{i=1}^w C_i\,(\ee_i,\ee_j)
		=C_j. \qedhere
		\]
	\end{proof}

	\section{Ultraweak closedness of the truncated \ratQM }
	\label{s:weighted-closed}

	We prove that the truncated rational quadratic
	modules are closed in the product ultraweak topology
	(Theorem~\ref{t:weighted-uw}).
	The key ingredient is a uniform bound on their
	Gram--CP representations
	(Lemma~\ref{l:weighted-bound}), obtained using
	the reproducing tuples constructed in
	Section~\ref{s:abhay-repro}.
	
	\subsection{Product ultraweak topology and intrinsic
		quadratic modules}
	Let $\cK$ be a complex Hilbert space. For $n\ge0$, let $\cA_n=\cB(\cK)\otimes V_n$,
	with the algebraic tensor product. Choose a basis $\bb b_1,\ldots,\bb b_s$ of $V_n$.
	A net $\sum_a C_{\alpha,a}\otimes\bb b_a$ converges in the
	\emph{product ultraweak topology} if every coefficient
	$C_{\alpha,a}$ converges ultraweakly in $\cB(\cK)$.
	This topology is independent of the choice of basis.
	Identifying $\cA_n$ with $\cB(\cK)^s$, it is the weak-$*$
	topology arising from the predual consisting of $s$-tuples of trace-class operators $\cT(\cK)^s$, with its
	$\ell^1$-norm. We write
	\[
	\left\|\sum_a C_a\otimes\bb b_a\right\|_{\rm c}
	:=\max_a\|C_a\|.
	\]

	Let $\cH,\cK$ be nonzero complex Hilbert spaces.
	The coefficients of $L$ act on $\cH$, while those of the rational
	function act on $\cK$. We keep the admissible set $R$ fixed.
	For a basis $\ee_1,\ldots,\ee_m$ of $V_\ell$ define the
	\emph{intrinsic \ratQM } $Q_\ell(\cK,L;\cR )$ by
	\[
	Q_\ell(\cK,L;\cR ) \,\coloneq\, \{ \mathbf e_{\ell}^{*} (G+\psi(L))\mathbf e_{\ell} \ \colon \ G\in\cB(\cK^m)_+,\,
	\psi\in\CP(\cS_L,\cB(\cK^m))  \},
	\]
	where $\mathbf e_{\ell} = {\rm col} (\ee_{1} I_{\cK}, \ldots, \ee_{m} I_{\cK})$.
	Thus $Q_\ell(\C^d,L;\cR )=Q_\ell^{(d)}(L;\cR )$ by
	Proposition~\ref{p:finite-gram}. A change of basis in $V_\ell$
	only changes the Gram--CP data by congruence, so the cone
	is independent of the basis. Also,
	$Q_\ell(\cK,L;\cR )\subseteq
	(\cB(\cK)\otimes V_{2\ell+1})_{\rm sa}$.
	
	\subsection{Uniform bounds and ultraweak closedness}
	\label{s:weighted-bounds}
	
	\begin{lem}\label{l:sep-eval}
		Let $k\ge0$ and let $\cW\subseteq V_k$ have basis
		$\ee_1,\ldots,\ee_w$. There are
		$n\ge1$, a hermitian tuple $X\in\bigcap_{r\in R}\hdom_n r$,
		and vectors $\eta,\zeta_1,\ldots,\zeta_w\in\C^n$ such that
		\[
		L(X)\succeq\tfrac12 I,\qquad
		\langle\ee_i(X)\eta,\zeta_j\rangle=\delta_{ij}.
		\]
		In particular, the functionals
		$\Omega_j(T)=\langle T\eta,\zeta_j\rangle$
		extract the coefficients of every element of
		$\cB(\cK)\otimes \cW$, for any Hilbert space $\cK$.
	\end{lem}
	\begin{proof}
		Let $\mathfrak A = \mathfrak A^{(\max\{1,k\})}$ be the tuple furnished by
		Proposition~\ref{p:reproducing-tuple} and put $\eta=\one$.
		The reproducing property shows that the vectors
		$\ee_1( \mathfrak A)\eta,\ldots,\ee_w(\mathfrak A)\eta$ are linearly
		independent. We claim that this remains true after replacing $\mathfrak A$
		by a sufficiently small positive scalar multiple, provided one
		avoids finitely many scalars.
		
		Indeed, along the line $t\mapsto t \mathfrak A$, each determinant needed to
		evaluate an expression in $R$ is a rational function of $t$.
		These determinants do not vanish at $t=1$, and hence have only
		finitely many zeros. The same reasoning applies to a nonzero
		minor of the matrix with columns $\ee_i(t \mathfrak A)\eta$.
		Thus the expressions in $R$ remain defined and the indicated
		vectors remain independent for arbitrarily small $t>0$.
		As $L$ is monic, we may also arrange that
		\[
		L(t \mathfrak A)=I+t\sum_{j=1}^{\gen}H_j\otimes \mathfrak A_j\succeq\tfrac12 I.
		\]
		Set $ X = t \mathfrak A$ and choose $\zeta_j$ dual to the independent
		vectors $\ee_i(X)\eta$. For
		$\mF=\sum_i C_i\otimes\ee_i$, we then have
		\[
		(\id\otimes\Omega_j)\mF(X)
		=\sum_i C_i\langle\ee_i(X)\eta,\zeta_j\rangle = C_j.
		\qedhere
		\]
	\end{proof}
	
	\begin{lem}\label{l:weighted-bound}
		Fix $\ell\ge0$ and a basis $\ee_1,\ldots,\ee_m$ of $V_\ell$.
		There is a constant $\kappa >0$, depending only on $L$, $R$, and the
		chosen bases, such that every representation
		\[
		\mF=\mathbf e_\ell^*(G+\psi(L))\mathbf e_\ell,
		\qquad G\succeq0,\quad\psi\in\CP(\cS_L,\cB(\cK^m)),
		\]
		satisfies
		\begin{equation}\label{e:weighted-bound}
			\|G\|+\|\psi(I_{\cH})\|\le \kappa \|\mF\|_{\mathrm c}.
		\end{equation}
		Here the coefficient norm is taken in a fixed basis of
		$V_{2\ell+1}$. In particular, $\kappa$ is independent of $\cK$.
	\end{lem}
	\begin{proof}
		Choose $X,\eta$ from Lemma~\ref{l:sep-eval} with $\cW=V_\ell$.
		The vectors $\ee_i(X)\eta$ are linearly independent. Thus
		there exists $\beta >0$ such that
		\begin{equation}\label{e:lower-extract}
			\Big\|\sum_i v_i\otimes\ee_i(X)\eta\Big\|^2
			\ge \beta \sum_i\|v_i\|^2
		\end{equation}
		for vectors $v_i$ in any Hilbert space. The same $\beta$ works
		for every such space.
		
		Factor $G=\mathsf R^* \mathsf R$, with block entries $\mathsf R_{ji}$. By Arveson's
		extension theorem and Stinespring's theorem, write
		$\psi(B) = \Gamma^*\rho(B) \Gamma$ for $B\in\cS_L$, where $\rho$ is a
		unital representation of $\cB(\cH)$ and
		$\Gamma = (\Gamma_1\ \cdots\ \Gamma_m)$. Set
		$\widetilde{\qq}=\sum_i \Gamma_i\otimes\ee_i$.
		Then
		\[
		\mathbf e_\ell^*\psi(L)\mathbf e_\ell
		=\widetilde{\qq}^{\,*}\rho(L)\widetilde{\qq}.
		\]
		As $L(X)\succeq\frac12I$, positivity is preserved by
		$\rho\otimes\id$, and hence
		$\rho(L)(X)\succeq\frac12I$. For $\xi\in\cK$,
		\begin{align*}
			&\langle(\mathbf e_\ell^*G\mathbf e_\ell)(X)
			(\xi\otimes\eta),\xi\otimes\eta\rangle\\
			&\qquad=\sum_j\Big\|\sum_i \mathsf R_{ji}\xi\otimes
			\ee_i(X)\eta\Big\|^2
			\ge \beta \sum_{i,j}\| \mathsf R_{ji}\xi\|^2
			=\beta \Big\langle\sum_iG_{ii}\xi,\xi\Big\rangle,
		\end{align*}
		while
		\begin{align*}
			&\langle(\mathbf e_\ell^*\psi(L)\mathbf e_\ell)(X)
			(\xi\otimes\eta),\xi\otimes\eta\rangle\\
			&\qquad\ge\frac12\Big\|\sum_i \Gamma_i\xi\otimes
			\ee_i(X)\eta\Big\|^2
			\ge\frac{\beta}{2}\sum_i\|\Gamma_i\xi\|^2
			=\frac{\beta}{2}\Big\langle
			\sum_i\psi(I_{\cH})_{ii}\xi,\xi\Big\rangle.
		\end{align*}
		Adding these estimates yields
		\[
		\beta \Big\langle
		\Big(\sum_iG_{ii}+\tfrac12\sum_i\psi(I_{\cH})_{ii}\Big)
		\xi,\xi\Big\rangle
		\le\|\eta\|^2\|\mF(X)\|\|\xi\|^2.
		\]
		For a positive block matrix, its norm is at most the sum
		of the norms of its diagonal blocks. We therefore obtain
		$\|G\|+\|\psi(I_{\cH})\|\le \kappa_1\|\mF(X)\|$.
		Since $X$ is fixed, evaluation at $X$ is bounded in the
		coefficient norm, independently of $\cK$. Thus
		$\|\mF(X)\|\le \kappa_2\|\mF\|_{\mathrm c}$,
		which proves the assertion.
	\end{proof}
	
	\begin{thm}\label{t:weighted-uw}
		For every $\ell\ge0$ and every Hilbert space $\cK$,
		$Q_\ell(\cK,L;\cR )$ is closed in the product ultraweak
		topology of $(\cB(\cK)\otimes V_{2\ell+1})_{\rm sa}$.
	\end{thm}
	\begin{proof}
		Choose a basis $\bb b_1,\ldots,\bb b_s$ of $V_{2\ell+1}$.
		It suffices, by the Krein--\v{S}mulian theorem, to show
		that the intersection of the cone with every closed ball
		of the coefficient norm is product-ultraweakly closed.
		
		Let $(\mF_\alpha)$ be a coefficient-norm bounded net in
		$Q_\ell(\cK,L;\cR )$ converging product-ultraweakly to
		$\mF$. Choose Gram--CP data $(G_\alpha,\psi_\alpha)$ for
		each $\mF_\alpha$. Lemma~\ref{l:weighted-bound} gives
		uniform bounds on $G_\alpha$ and $\psi_\alpha(I_{\cH})$.
		Since $\psi_\alpha$ is completely positive,
		$\|\psi_\alpha\|=\|\psi_\alpha(I_{\cH})\|$.
		
		Let $B_1,\ldots,B_p$ be a basis of $\cS_L$.
		The Banach--Alaoglu theorem and the finiteness of this basis give a
		subnet, still indexed by $\alpha$, for which
		\[
		G_\alpha\longrightarrow G,\qquad
		\psi_\alpha(B_i)\longrightarrow \Psi_i\quad(1\le i\le p)
		\]
		ultraweakly. The operators $\Psi_i$ determine a linear map
		$\psi:\cS_L\to\cB(\cK^m)$; by linearity,
		$\psi_\alpha(B)\to\psi(B)$ ultraweakly for all
		$B\in\cS_L$. The ultraweak closedness of the positive cone
		gives $G\succeq0$. Applied to positive matrices over
		$\cS_L$, the same argument gives
		$[\psi(B_{ij})]\succeq0$ whenever
		$[B_{ij}]\succeq0$. Thus $\psi$ is completely positive.
		
		To identify the limit, expand
		\[
		\mF_\alpha
		=\sum_{i,j}\bigl((G_\alpha)_{ij}
		+\psi_\alpha(I_{\cH})_{ij}\bigr)
		\otimes\ee_i^*\ee_j
		+\sum_{k=1}^{\gen}\sum_{i,j}
		\psi_\alpha(H_k)_{ij}\otimes\ee_i^*x_k\ee_j.
		\]
		Every rational function on the right belongs to
		$V_{2\ell+1}$ and is independent of $\alpha$.
		Upon expressing these functions in the fixed basis
		$\bb b_1,\ldots,\bb b_s$, each coefficient is a finite
		linear combination of blocks of $G_\alpha$ and of
		$\psi_\alpha(I_{\cH}),\psi_\alpha(H_1),\ldots,
		\psi_\alpha(H_{\gen})$. Their ultraweak limits are
		therefore the coefficients of
		$\mathbf e_\ell^*(G+\psi(L))\mathbf e_\ell$.
		By uniqueness of the product-ultraweak limit,
		\[
		\mF=\mathbf e_\ell^*(G+\psi(L))\mathbf e_\ell.
		\]
		Hence $\mF$ belongs to the cone. The assertion now
		follows from the Krein--\v{S}mulian theorem, since the cone is convex
		by Proposition~\ref{p:finite-gram} and its intrinsic
		Gram--CP description.
	\end{proof}
	
	\section{Finite-dimensional rational Positivstellensatz}\label{s:finite-fixed}

	In this section we prove the matrix-valued rational Positivstellensatz
	with an explicit truncation bound
	(Theorem~\ref{t:finite-fixed}).
	The proof combines separation of the truncated
	quadratic module with the construction of feasible
	finite-dimensional extensions of multiplication
	operators (Lemma~\ref{l:finite-feasible-extension}).

	\begin{lem}\label{l:finite-gns-positive}
		Let $\ell\ge2$, and let
		$\lambda: \opm_d(V_{2\ell+1})   \to \C$ be a linear
		functional which is strictly positive on
		$Q_\ell^{(d)}(L;\cR )\setminus\{0\}$.
		Then
		\[
		\langle\uu,\rv\rangle=\lambda(\rv^*\uu)
		\]
		defines an inner product on
		$U=V_\ell^{1\times d}$.
		The multiplication operators
		$\fX_j:S=V_{\ell-1}^{1\times d}\to U$
		are symmetric on $S$. If
		$A_j=P_S\fX_j|_S$, then
		\[
		L( A)=I_{\cH}\otimes I_S+
		\sum_{j=1}^\gen H_j\otimes  A_j \succeq \delta I
		\]
		for some $\delta >0.$
	\end{lem}
	
	\begin{proof}
		If $\uu\ne0$, then $\uu^*\uu$ is a nonzero
		element of $Q_\ell^{(d)}(L;\cR )$, so
		$\lambda(\uu^*\uu)>0$.
		Moreover, for $\uu,\rv\in S$,
		\[
		\langle\fX_j\uu,\rv\rangle
		=\lambda(\rv^*x_j\uu)
		=\langle\uu,\fX_j\rv\rangle.
		\]
		
		Given $0\ne\zeta\in\cH\otimes S$, write
		\[
		\zeta=\sum_{a=1}^k h_a\otimes\uu_a,
		\]
		where the $h_a$ are orthonormal and
		$\uu_a\in S$. We may assume $k < m d$, where $m =  \dim V_\ell$.
		Let $\Theta:\C^k\to\cH$ be the isometry defined by
		$\Theta e_a=h_a$. Fix a state $\omega$ on $\cS_L$ and define
		\[
		\pi(B)=\Theta^* B \Theta\oplus\omega(B)I_{md-k},
		\qquad B\in\cS_L.
		\]
		Then $\pi\in\UCP(\cS_L,\opm_{md}(\C))$.
		Let $\qq\in\opm_{m d\times d}(V_{\ell-1})$
		have rows $\uu_1,\ldots,\uu_k$, followed
		by zero rows. Then
		\begin{align*}
			\langle L( A)\zeta,\zeta\rangle
			&=
			\sum_{a=1}^k
			\langle\uu_a,\uu_a\rangle
			+
			\sum_{j=1}^{\gen}\sum_{a,b=1}^k
			\langle H_jh_b,h_a\rangle
			\langle  A_j\uu_b,\uu_a\rangle\\
			&=
			\sum_{a=1}^k
			\lambda(\uu_a^*\uu_a)
			+
			\sum_{j=1}^{\gen}\sum_{a,b=1}^k
			\langle H_jh_b,h_a\rangle
			\lambda(\uu_a^*x_j\uu_b)\\
			&=
			\lambda\left(
			\sum_{a,b=1}^k
			\pi(I_{\cH})_{ab}\uu_a^*\uu_b
			+
			\sum_{j=1}^{\gen}\sum_{a,b=1}^k
			\pi(H_j)_{ab}\uu_a^*x_j\uu_b
			\right)\\
			&=
			\lambda\bigl(\qq^*\pi(L)\qq\bigr)>0.
		\end{align*}
		The inequality is strict because
		$\qq^*\pi(L)\qq$ is a nonzero element of
		$Q_\ell^{(d)}(L;\cR )$. Indeed, since $\zeta\ne0$,
		we have $\qq\ne0$. By Lemma~\ref{l:sep-eval},
		there exists a hermitian tuple $X\in\bigcap_{r\in R}\hdom_n r$  such that $L(X)\succeq\frac12 I$
		and evaluation at $X$ is injective on $V_\ell$.
		In particular, $\qq(X)\ne0$. Since $\pi$ is UCP,
		$\pi(L)(X)\succeq\frac12 I$, and hence
		\[
		(\qq^*\pi(L)\qq)(X)
		=\qq(X)^*\pi(L)(X)\qq(X)
		\succeq\frac12\qq(X)^*\qq(X)\ne0.
		\]
		Thus $\qq^*\pi(L)\qq\ne0$, and the strict
		positivity of $\lambda$ on
		$Q_\ell^{(d)}(L;\cR )\setminus\{0\}$ proves the claim.
		
		Since the unit sphere of
		$\opm_{md \times d}(V_{\ell-1})$ and
		$\UCP(\cS_L,\opm_{md}(\C))$ are compact,
		the strict positivity above is uniform
		for $\|\zeta\|=1$.
		Thus $L(  A)\succeq\delta I$ for some
		$\delta>0$.
	\end{proof}

	\begin{lem}\label{l:finite-feasible-extension}
		Let $\ell\ge2$, and equip
		$U=V_\ell^{1\times d}$ with an inner product
		for which the compressions
		\[
		A_j=P_S\fX_j|_S,
		\qquad S=V_{\ell-1}^{1\times d},
		\]
		are hermitian. Suppose that $L(A)\succeq \delta I$ for some $\delta > 0$.
		Then there exists a finite-dimensional
		hermitian tuple
		\[
		X\in\cD(L)\cap
		\bigcap_{r\in R}\hdom r
		\]
		such that
		\begin{equation}\label{e:finite-extension-agreement}
			X_j|_{V_{\ell-2}^{1\times d}}
			=\fX_j|_{V_{\ell-2}^{1\times d}}
			\qquad(1\le j\le\gen).
		\end{equation}
	\end{lem}
	
	\begin{proof}
		We first work with the inner product induced by $\phi$ in
		Section~\ref{s:mult}. By Lemma~\ref{l:fX}, the
		multiplication operators are hermitian and lie in
		$\bigcap_{r\in R}\dom_\infty r$.

		Write $U=S\oplus_\phi U_0$, where
		$U_0$ is the orthogonal complement of $S$ in $U$
		with respect to the inner product induced by $\phi$,
		and let $\kappa =\dim U_0$.
		By Proposition~\ref{p:main}, there exist $n\ge \kappa$,
		an invertible $E_0\in\mat{n}$, and a hermitian tuple
		\begin{equation} \label{eq:matrix with respect to phi}
			\widehat X_j=
			\begin{pmatrix}
				A_j^\phi & (Y_j^*\ \ 0)E_0^*\\
				E_0\binom{Y_j}{0} & \widetilde Z_j
			\end{pmatrix},
			\qquad 1\le j\le\gen,
		\end{equation}
		on $S\oplus\C^n$ such that
		$\widehat X\in\bigcap_{r\in R}\hdom r$.
		Here
		\[
		\fX_j|_S=
		\binom{A_j^\phi}{Y_j}
		\]
		is written relative to $U=S\oplus_\phi U_0$, and $A_j^\phi$ denotes the compression with respect to the $\phi$-inner product.
		
		Identifying $\C^n$ with $U_0\oplus\C^{n-\kappa}$,
		define
		\[
		J:U\longrightarrow S\oplus\C^n,\qquad
		J(s+u)=s\oplus E_0\binom{u}{0},
		\quad s\in S,\quad u\in U_0.
		\]
		Since $E_0$ is invertible, $J$ is injective.
		Moreover, for $s\in S$,
		\[
		\widehat X_jJs
		=
		\binom{A_j^\phi s}{E_0\binom{Y_js}{0}}
		=
		J\fX_js.
		\]
		Thus,
		\[
		\widehat X_jJ|_S=J\fX_j|_S,
		\qquad 1\le j\le\gen.
		\]
		
		Choose an orthonormal basis of $S$ and extend it
		to an orthonormal basis of $U$, with respect to
		the given inner product (not the one  induced by $\phi$). Applying $J$ to these basis
		vectors and extending them gives a basis of
		$S\oplus\C^n$, which is not required to be orthonormal.
		The matrices $T_j$ of $\widehat X_j$ in this
		basis are simultaneously similar to
		$\widehat X_j$ in \eqref{eq:matrix with respect to phi} and have the form
		\[
		T_j=
		\begin{pmatrix}
			A_j&*&*\\
			Y_j&*&*\\
			0&*&*
		\end{pmatrix},
		\qquad A_j=A_j^*,
		\]
		where
		\[
		\fX_j|_S=
		\begin{pmatrix}A_j\\Y_j\end{pmatrix}
		\]
		is written relative to the orthogonal
		decomposition $U=S\oplus S^\perp$
		for the given inner product.

		By Lemma~\ref{l:change-basis-completion},
		there exist an invertible matrix $E$ and a
		hermitian tuple $Z$ such that
		\[
		\widetilde X_j=
		\begin{pmatrix}
			A_j&(Y_j^*\ \ 0)E^*\\
			E\binom{Y_j}{0}&Z_j
		\end{pmatrix},
		\qquad 1\le j\le\gen,
		\]
		belongs to $\bigcap_{r\in R}\hdom r$.

		For $\varepsilon_0,\varepsilon_1>0$, let
		\[
		X_j(\varepsilon_0,\varepsilon_1)=
		\begin{pmatrix}
			A_j&
			\varepsilon_0(Y_j^*\ \ 0)E^*\\
			\varepsilon_0E\binom{Y_j}{0}&
			\varepsilon_1Z_j
		\end{pmatrix}.
		\]
		Since $X(1,1)=\widetilde X$ belongs to
		$\bigcap_{r\in R}\hdom r$, the product of the
		denominator determinants of the expressions in $R$,
		evaluated at $X(\varepsilon_0,\varepsilon_1)$,
		is a polynomial that does not vanish at $(1,1)$.
		Thus, it is not identically zero and cannot
		vanish  entirely on any open subset of $\R^2$.
		Hence, for arbitrarily small positive
		$\varepsilon_0,\varepsilon_1$, we have
		\[
		X(\varepsilon_0,\varepsilon_1)
		\in\bigcap_{r\in R}\hdom r.
		\]
		
		Since $L(A)\succeq \delta I$ for some $\delta > 0$, we have
		\[
		L(X(0,0))=L(A)\oplus I\succeq 
		\min\{\delta,1\} I. 
		\]
		By continuity, $L(X(\varepsilon_0,\varepsilon_1))$
		remains positive definite for sufficiently small
		$\varepsilon_0,\varepsilon_1>0$. Thus, 
		we can choose these parameters so that
		\[
		X:=X(\varepsilon_0,\varepsilon_1)
		\in\cD(L)\cap\bigcap_{r\in R}\hdom r.
		\]
		
		Finally,
		$\fX_j(V_{\ell-2}^{1\times d})\subseteq S$,
		so $Y_j$ vanishes on $V_{\ell-2}^{1\times d}$. This proves
		\eqref{e:finite-extension-agreement}.
	\end{proof}
	
	\begin{thm}\label{t:finite-fixed}
		Let
		$ \mF=\sum_{a=1}^{N}C_a\otimes\rr_a
		\in\opm_d(\C)\otimes V_1 $
		be hermitian, where $C_a\in\opm_d(\C)$ and
		$r_a\in R$ represents $\rr_a$ for each $a$.
		Set
		$	\ell_*=2\tau(R)+1 $.
		Then $\mF(X)\succeq0$ for every
		$ X\in\cD(L)\cap\bigcap_{r\in R}\hdom r $
		if and only if
		$	\mF\in Q_{\ell_*}^{(d)}(L;\cR ). $
	\end{thm}
	
	\begin{proof}
		The reverse implication follows from complete positivity.
		For the forward implication, suppose that
		$\mF\notin Q_{\ell_{*}}^{(d)}(L;\cR )$.
		By Theorem~\ref{t:weighted-uw}, the cone
		$Q_{\ell_*}^{(d)}(L;\cR )$ is closed. It is also pointed.
		Indeed, Lemma~\ref{l:sep-eval}, applied to
		$V_{2\ell_{*}+1}$, gives a hermitian tuple $X$ such that
		$L(X)\succ0$ and evaluation at $X$ is injective on
		$V_{2\ell_*+1}$. If both $\mG$ and $-\mG$ belong to
		$Q_{\ell_*}^{(d)}(L;\cR )$, then
		\[
		\mG(X)\succeq0
		\quad\text{and}\quad
		-\mG(X)\succeq0.
		\]
		Hence $\mG(X)=0$, and injectivity gives $\mG=0$.

		Since $Q_{\ell_*}^{(d)}(L;\cR )$ is a closed pointed
		convex cone in a finite-dimensional real vector
		space, the Hahn--Banach separation theorem
		(see \cite[Theorem III.1.3]{Bar}) gives a real
		linear functional
		$	\lambda:\opm_d(V_{2\ell_*+1})_{\rm sa}\to\R $
		such that
		\[
		\lambda(\mF)<0,
		\qquad
		\lambda(\mG)>0
		\quad
		(0\ne\mG\in Q_{\ell_*}^{(d)}(L;\cR )).
		\]
		We extend $\lambda$ complex linearly to
		$\lambda:\opm_d(V_{2\ell_*+1})\to\C$.

		Equip $U=V_{\ell_*}^{1\times d}$ with the inner product
		defined in Lemma~\ref{l:finite-gns-positive}.
		By that lemma, the compressed multiplication tuple
		$A$ is hermitian and satisfies $L(A)\succeq\delta I$
		for some $\delta>0$.
		By Lemma~\ref{l:finite-feasible-extension},
		we obtain a finite-dimensional hermitian tuple
		$	X\in\cD(L)\cap\bigcap_{r\in R}\hdom r $
		such that
		\[
		X_j|_{V_{\ell_*-2}^{1\times d}}
		=\fX_j|_{V_{\ell_*-2}^{1\times d}},
		\qquad j=1,\ldots,\gen.
		\]

		Since
		\[
		2\tau(r_a)\le 2\tau(R)\le\ell_*=(\ell_*-2)+2,
		\]
		Proposition~\ref{p:gns}, applied with
		$\ell=\ell_*-2$ and $s=1$, gives
		\[
		r_a(X)\rv_i=\rr_a\rv_i,
		\qquad i=1,\ldots,d,
		\]
		where $\rv_i\in V_0^{1\times d}$ is the row
		vector with $1$ in the $i$-th position and
		zeros elsewhere.

		Let $e_1,\ldots,e_d$ be the standard basis of $\C^d$
		and set
		$	\xi=\sum_{i=1}^d e_i\otimes\rv_i. $
		Using the preceding identities and the definition
		of the inner product, we obtain
		\begin{align*}
			\langle\mF(X)\xi,\xi\rangle
			&=\sum_{a=1}^N\sum_{i,j=1}^d
			(C_a)_{ij}\langle r_a(X)\rv_j,\rv_i\rangle\\
			&=\sum_{a=1}^N\sum_{i,j=1}^d
			(C_a)_{ij}\lambda(\rv_i^*\rr_a\rv_j)\\
			&=\lambda\left(
			\sum_{a=1}^N C_a\otimes\rr_a
			\right)\\
			&=\lambda(\mF)<0.
		\end{align*}
		This contradicts the assumed positivity of $\mF(X)$.
		Therefore $\mF\in Q_{\ell_*}^{(d)}(L;\cR )$.
	\end{proof}
	
	\section{From finite-dimensional to operator-valued rational functions}
	\label{s:infinite-K}

	We extend the finite-dimensional Positivstellensatz
	(Theorem~\ref{t:finite-fixed})
	to rational functions with coefficients in
	$\cB(\cK)$ (Theorem~\ref{t:operator-K}).
	The passage to infinite-dimensional coefficient spaces
	relies on the uniform truncation bound of
	Theorem~\ref{t:finite-fixed} and the ultraweak closedness
	established in Theorem~\ref{t:weighted-uw}.

	\begin{thm}\label{t:operator-K}
		Let $\cK$ be an infinite-dimensional Hilbert space,
		and let $\mF=\sum_{a=1}^N C_a\otimes\rr_a
		\in\cB(\cK)\otimes V_1$ be hermitian, where
		$C_a\in\cB(\cK)$ and $r_a\in R$ represents $\rr_a$.
		Set $\ell_*=2\tau(R)+1$.
		If $\mF(X)\succeq0$ for every
		$X\in\cD(L)\cap\bigcap_{r\in R}\hdom r$,
		then $\mF\in Q_{\ell_*}(\cK,L;\cR )$.
		Equivalently, there exist
		$\rr,\qq\in\cB(\cK)\otimes V_{\ell_*}$ and
		$\pi\in\UCP(\cS_L,\cB(\cK))$ such that
		\begin{equation}\label{e:single-infinite}
			\mF=\rr^*\rr+\qq^*\pi(L)\qq.
		\end{equation}
		The converse also holds.
	\end{thm}
	\begin{proof}
		Let $(P_\alpha)$ be the net of finite-rank orthogonal
		projections on $\cK$, ordered by inclusion of their ranges.
		Set $\cK_\alpha=P_\alpha\cK$ and
		$d_\alpha=\dim\cK_\alpha$, and consider
		\[
		\mF_\alpha
		=\sum_{a=1}^N
		(P_\alpha C_aP_\alpha)|_{\cK_\alpha}\otimes\rr_a.
		\]
		Since compression preserves positivity,
		$\mF_\alpha(X)\succeq0$ for every
		$X\in\cD(L)\cap\bigcap_{r\in R}\hdom r$.
		The expressions $r_a$ and the admissible set $R$
		remain unchanged, so the bound $\ell_*$ is the same
		for every $\alpha$. By Theorem~\ref{t:finite-fixed},
		$	\mF_\alpha\in Q_{\ell_*}^{(d_\alpha)}(L;\cR ) $.

		{Let $m = \dim V_{\ell_{*}}$}
		and let $\iota_\alpha:\cK_\alpha^m\to\cK^m$ denote the
		natural inclusion. Choose Gram--CP data
		$(G_\alpha,\psi_\alpha)$ for $\mF_\alpha$ and define
		\[
		\widetilde G_\alpha=\iota_\alpha G_\alpha\iota_\alpha^*,
		\qquad
		\widetilde\psi_\alpha(B)
		=\iota_\alpha\psi_\alpha(B)\iota_\alpha^*,
		\quad B\in\cS_L.
		\]
		Then $\widetilde G_\alpha\succeq0$ and
		$\widetilde\psi_\alpha$ is completely positive. Moreover,
		\[
		\sum_{a=1}^N P_\alpha C_aP_\alpha\otimes\rr_a
		=
		\mathbf e_{\ell_*}^*
		\bigl(\widetilde G_\alpha+
		\widetilde\psi_\alpha(L)\bigr)
		\mathbf e_{\ell_*}
		\in Q_{\ell_*}(\cK,L;\cR ).
		\]
		Since $P_\alpha\to I_{\cK}$ strongly,
		$P_\alpha C_aP_\alpha\to C_a$ ultraweakly for
		each $a$. Thus the expression on the left converges
		to $\mF$ in the product ultraweak topology.
		By Theorem~\ref{t:weighted-uw},
		\[
		\mF\in Q_{\ell_*}(\cK,L;\cR ).
		\]

		Write
		$	\mF=\mathbf e_{\ell_*}^*(G+\psi(L))\mathbf e_{\ell_*}, $
		where $G\succeq0$ and $\psi$ is completely positive.
		Set $\mathsf R=G^{1/2}$ and
		$D=\psi(I_{\cH})^{1/2}$.
		As in the proof of Proposition~\ref{p:finite-gram},
		there is a UCP map
		$\widehat\pi:\cS_L\to\cB(\cK^m)$ such that
		$\psi(B)=D\widehat\pi(B)D$ for every $B\in\cS_L$.
		
		Since $\cK$ is infinite-dimensional, we may choose
		a unitary $J:\cK^m\to\cK$.
		Using the basis of $V_{\ell_*}$, the block entries
		of $J\mathsf R$ and $JD$ determine elements
		$\rr,\qq\in\cB(\cK)\otimes V_{\ell_*}$.
		Taking $\pi(B)=J\widehat\pi(B)J^*$, we obtain
		\[
		\mF=\rr^*\rr+\qq^*\pi(L)\qq,
		\]
		as required.
	\end{proof}
	
	\subsection{Operator-valued noncommutative rational Hilbert's 17th problem}

	Taking $L=I$ in the preceding Positivstellens\"atze
	yields operator-valued versions of noncommutative
	rational Hilbert's 17th problem.
	We obtain a single-square representation for
	infinite-dimensional coefficient spaces
	(Corollary~\ref{c:operator-H17}) and a finite
	sum-of-squares representation in the matrix-valued
	case (Corollary~\ref{c:finite-H17}).

	\begin{cor}
		\label{c:operator-H17}
		Let $\cK$ be an infinite-dimensional Hilbert space and let
		\[
		\mF=\sum_{a=1}^N C_a\otimes\rr_a
		\in\cB(\cK)\otimes\rx
		\]
		be hermitian, where $r_a\in\rex$ represents $\rr_a$.
		If $\mF(X)\succeq0$ for every
		\[
		X\in\bigcap_{a=1}^N\hdom r_a,
		\]
		then there exists $\rs\in\cB(\cK)\otimes\rx$ such that
		\[
		\mF=\rs^*\rs.
		\]
	\end{cor}
	
	\begin{proof}
		Choose a finite admissible set $R$ containing $r_1,\ldots,r_N$. 
		Applying Theorem~\ref{t:operator-K}
		with $L=I$ gives
		\[
		\mF=\rr^*\rr+\qq^*\qq
		\]
		for some $\rr,\qq\in\cB(\cK)\otimes\rx$.
		Since $\cK$ is infinite-dimensional, there exist
		isometries $U,V\in\cB(\cK)$ with $U^*V=0$.
		Thus $\rs=U\rr+V\qq$ satisfies $\mF=\rs^*\rs$.
	\end{proof}

	\begin{cor}
		\label{c:finite-H17}
		Let $\cK$ be a nonzero finite-dimensional Hilbert space
		and let
		\[
		\mF=\sum_{a=1}^N C_a\otimes\rr_a
		\in\cB(\cK)\otimes\rx
		\]
		be hermitian, where $r_a\in\rex$ represents $\rr_a$.
		If $\mF(X)\succeq0$ for every
		\[
		X\in\bigcap_{a=1}^N\hdom r_a,
		\]
		then there exist finitely many
		$\rs_j\in\cB(\cK)\otimes\rx$ such that
		\begin{equation}\label{eq:H17sos}
			\mF=\sum_j\rs_j^*\rs_j.
		\end{equation}
		Moreover, for any finite admissible set $R$
		containing $r_1,\ldots,r_N$, one may take
		$\rs_j\in\cB(\cK)\otimes V_{\ell_*}$,
		with at most $\dim V_{\ell_*}$ summands,
		where $\ell_*=2\tau(R)+1$.
	\end{cor}
	
	\begin{proof}
		Fix a finite admissible set $R$ containing
		$r_1,\ldots,r_N$.
		Identify $\cB(\cK)$ with $\opm_d(\C)$, where
		$d=\dim\cK$. Applying
		Theorem~\ref{t:finite-fixed} with $L=I$ gives
		$\mF\in Q_{\ell_*}^{(d)}(I;\cR)$.
		By Proposition~\ref{p:finite-gram},
		\[
		\mF=\mathbf e_{\ell_*}^*H\mathbf e_{\ell_*},
		\qquad H\succeq0.
		\]
		Writing $H=T^*T$ and taking the block rows
		$\rs_j$ of $T\mathbf e_{\ell_*}$ yields
		the required \eqref{eq:H17sos}.
	\end{proof}
	
	\subsection{The commutative case}\label{ssec:comm-case}
	
	The operator-valued analog of Hilbert's 17th problem
	fails in the commutative setting, even for polynomials.
	
	Indeed, let $\mathcal P$ be the cone of nonnegative
	ternary sextic forms, and choose a sequence $(p_n)$
	dense in its unit-norm section. Then
	\[
	F=\operatorname{diag}(p_1,p_2,\ldots)
	\in \cB(\ell^2)\otimes\C[x_1,x_2,x_3]
	\]
	is positive semidefinite everywhere. If
	$F=\sum_j q_j^*q_j$ for finitely many
	$q_j\in\cB(\ell^2)\otimes\C(x_1,x_2,x_3)$,
	clearing their finitely many scalar denominators
	gives a nonzero $h\in\R[x_1,x_2,x_3]$ such that
	\[
	h^2p_n\in\sum\R[x_1,x_2,x_3]^2
	\qquad(n\ge1).
	\]
	By density and closedness of the sum-of-squares cone
	in fixed degree, the same holds for every
	$p\in\mathcal P$. But this
	contradicts Reznick's theorem on the absence
	of uniform denominators \cite{Rez05}.
	
	In contrast, the matrix-valued commutative
	Hilbert's 17th problem has an affirmative answer.
	For every $d$, a hermitian matrix
	$F\in \opm_d(\C(x_1,\ldots,x_{\gen}))$ that is positive
	semidefinite wherever defined on $\R^{\gen}$
	is a finite sum of hermitian squares of rational
	matrices. The real symmetric version was established by
	Gondard and Ribenboim \cite{GR74} and independently
	by Procesi and Schacher \cite{PS76}.
	The complex hermitian version is due to
	Djokovi\'c \cite{Djo76}; see also the
	constructive proof of Hillar and Nie \cite{HN08}.
	Indeed, rational congruence diagonalization reduces
	the assertion to Artin's scalar theorem, applied to
	the nonnegative diagonal entries.

\end{document}